# HYPERBOLIC 3-MANIFOLDS OF BOUNDED VOLUME AND TRACE FIELD DEGREE

BOGWANG JEON

ABSTRACT. For a single cusped hyperbolic 3-manifold, Hodgson proved that there are only finitely many Dehn fillings of it whose trace fields have bounded degree. In this paper, we conjecture the same for manifolds with more cusps, and give the first positive results in this direction. For example, in the 2-cusped case, if a manifold has linearly independent cusp shapes, we show that the manifold has the desired property. To prove the results, we use Habegger's proof of the Bounded Height Conjecture in arithmetic geometry.

## 1. Introduction

In the study of hyperbolic 3-manifolds, the following question is very natural.

**Question 1.** For a given number $D > 0$, are there only finitely many hyperbolic 3-manifolds whose volumes and degrees of their trace fields are bounded by $D$?

By the Jorgensen-Thurston theory (see Theorem 2.5 in Section 2.3), to answer **Question 1**, it is enough to answer the following question:

**Question 2.** For a $k$-cusped manifold $M$ and a constant $D > 0$, are there only finitely many Dehn fillings of $M$ whose trace fields have degree $\leq D$?

It is commonly believed that the answer to both questions is yes and this was proved for the 1-cusped case by Hodgson (see [11] for a relevant more generalized version), but little was known for manifolds with $k \geq 2$ cusps in general. In this paper we answer these questions for special types of manifolds with more cusps. For instance, the following is one of the main theorems of this paper.

**Theorem 1.1.** *Let $M$ be a 2-cusped hyperbolic 3-manifold having cusp shapes $\tau_1$ and $\tau_2$. If $1, \tau_1, \tau_2, \tau_1\tau_2$ are linearly independent over $\mathbb{Q}$, then, for any $D > 0$, there are only finitely many Dehn fillings of $M$ whose trace field has degree less than $D$.*

From now on, for simplicity, we say $M$ has *rationally independent cusp shapes* if it satisfies the given condition in the above statement. Note that the linear independency of $1, \tau_1, \tau_2, \tau_1\tau_2$ over $\mathbb{Q}$ is independent of the choice of basis.

To prove the theorem, we first employ the notion of height from number theory, which is the standard way of measuring the complexity of algebraic numbers, and define it for each Dehn filling of $M$. Specifically, we define it as the trace value of the core geodesic of a Dehn filling. It is a fundamental theorem in number theory that there are only finitely many algebraic numbers of bounded height and degree. Hence, in terms of height instead of degree, to get the affirmative answer to **Question 2**, it is enough to answer the following stronger question (we'll discuss this in more detail in Section 3.2):

**Question 3.** For a $k$-cusped manifold $M$, is there a constant $D > 0$ such that, for any Dehn filling of $M$, its height is uniformly bounded by $D$?

According to Thurston's hyperbolic Dehn filling theory, each Dehn filled manifold of $M$ corresponds to a point on the deformation variety (of hyperbolic structures on $M$) satisfying certain additional conditions regarding to its Dehn filling coefficients. (For the moment, let's call this point on the deformation variety a "Dehn filling point". We'll give the precise definition later in Section 2.3.) By using the appropriate version of the deformation variety (precisely, the one having the holonomies of the longitude-meridian pairs as parameters), these conditions can be represented by a set of multiple equations defining an algebraic subgroup. So a Dehn filling point on the deformation variety becomes an intersection point between the deformation variety and an algebraic subgroup. Furthermore, using some elementary properties of height, it can be shown that if the height of a Dehn filling point is bounded, then the height of the corresponding Dehn filled manifold is also bounded. Thus, to answer **Question 3**, it is sufficient to prove the heights of intersection points (i.e. Dehn filling points) between the given algebraic varieties are uniformly bounded. As a result, the original problem in hyperbolic geometry is transformed into a problem in arithmetic geometry.

The height distribution of points on an algebraic variety is widely studied topic in arithmetic geometry and there are various theorems regarding to this theme. Among them, we use the one which is so called the Bounded Height Conjecture, originally formulated by E. Bombieri, D. Masser, U. Zannier in [3], and proved by P. Habegger in [5] (see also [13]).

**Theorem 1.2.** (Bounded Height Conjecture=Habegger's theorem) *Let $X \subset (\overline{\mathbb{Q}}^*)^n$ be an irreducible variety over $\overline{\mathbb{Q}}$. Then there is a Zariski open subset $X^{oa}$ of $X$, which is the complement of the union of anomalous subvarieties of $X$, so that the height is bounded in the intersection of $X^{oa}$ with the union of algebraic subgroups of dimension $\leq n - \dim X$.*

Since it takes quite a bit of background to define an anomalous variety, we postpone it until Section 3.3.

Here the point is that $X^{oa}$ is a Zariski open subset of $X$. Therefore, by applying the above theorem, we immediately get the uniformly boundedness of the height on most of $X$ unless $X^{oa} = \emptyset$ (but this can happen, unfortunately). However, there is a technical issue which prevents us from directly applying Habegger's theorem. Whereas we only interested in a local neighborhood of a point on the deformation variety, the bounded height conjecture deals with the whole variety. In addition, the singularity of the deformation variety may cause some unexpected problems which will be addressed in Section 5.1 deeply. Indeed, to get the desired result, we need to strengthen Theorem 1.2 a little. It turns out that, by following the original proof of Habegger's paper, we can extend the theorem in the way we can naturally apply to our situation (see Theorem 5.2).

Using this generalized version, we prove the following main theorem of the paper:

**Theorem 1.3.** *Suppose that the answer is yes to Question 3 for any $s$-cusped manifolds where $1 \leq s \leq k-1$. Let $X$ be the deformation variety of $k$-cusped hyperbolic 3-manifold $M$. If $X$ is simple, then the answer is yes to Question 3 for $M$.*

For the precise definition of a simple variety, see Definition 5.5 (the definitions in Section 3.3 are also needed). The definition of it is very natural. For instance, when $X$ is a 2-dimensional variety, it simply means $X^{oa}$ is nonempty, and, for the higher dimensional cases, the idea is extended in an analogous fashion.

Although we don't have any geometric criteria to judge when a hyperbolic 3-manifold has a simple deformation variety, we believe that "simple" is the general phenomenon. For example, if the deformation variety is not simple, we prove the following under the same assumption as in Theorem 1.1.

**Theorem 1.4.** *Let $M$ be a 2-cusped hyperbolic 3-manifold with rationally independent cusp shapes. If the deformation variety of $X$ is not simple, then the two cusps of $M$ are strongly geometrically isolated.*

Since strong geometric isolation is relatively rare, it is expected that the deformation variety being simple is quite common.

In the above case, if the two cusps of $M$ are strongly geometrically isolated, then $X^{oa} = \emptyset$ so we cannot apply the bounded height conjecture. However, in this case, interestingly enough, we can use Hodgson's method to show uniform boundedness of the heights. As a consequence, combining with Theorem 1.3, when a 2-cusped manifold has rationally independent cusp shapes, then whether its deformation variety is simple or not, the height of each Dehn filling is uniformly bounded (i.e. Theorem 1.1 holds).

For the higher cusped cases in general, the non-simple phenomenon is poorly understood, but we think Theorem 1.4 can be further extended, so we formulate it as a conjecture:

**Conjecture 1.** *Let $X$ be a $k$-cusped hyperbolic 3-manifold. If the deformation variety of $X$ is not simple, then $M$ has a set of cusps which are strongly geometrically isolated from the rest.*

This conjecture, together with Theorem 1.3, suggest the following seemingly plausible conjecture, which is the affirmative answer to **Question 3**:

**Conjecture 2.** (Bounded Height Conjecture in Hyperbolic 3-manifolds) *Let $M$ be a $k$-cusped hyperbolic 3-manifold. Then the height of any Dehn filling of $M$ is uniformly bounded.*

Even though we only deal with manifolds under certain restrictions, it is strongly believed that the above conjecture is true and this approach will eventually give us the complete positive answer to **Question 1**.

Lastly we exhibit an explicit example whose deformation variety is simple, but which is not covered by Theorem 1.1. Surely this also implies that most deformation varieties would be simple.

**Theorem 1.5.** *Let $W$ be the complement of the $(-2, 3, 8)$-pretzel link. Then the deformation variety of $W$ is simple.*

Here is the outline of the paper. In Section 2 and 3, we study some necessary background, and prove Theorem 1.4 and Theorem 1.3 in Sections 4 and 5 respectively. In Section 6, we show Theorem 1.5, and make some comments in Section 7. Finally we prove the generalized version of the Bounded Height Conjecture (Theorem 5.2) in Section 8.

## 2. Background I (Hyperbolic Geometry)

Before starting this section, let us note that we use the same notations repeatedly in different sections. That is, once we introduce a new notation, we will use it in later sections in the same meaning without defining it again.

**2.1. Gluing variety** In this section, we follow the same scheme in [9]. Suppose that $M$ is a $k$-cusped manifold whose hyperbolic structure is realized as a union of $n$ geometric tetrahedra having modulus $z_v$ $(1 \leq v \leq n)$. Then the gluing variety of $M$ is defined by the following form of $n$ equations where each represents the gluing condition at each edge of a tetrahedron:

$$(2.1) \qquad \prod_{v=1}^{n} z_v^{\theta_1(r,v)} \cdot (1 - z_v)^{\theta_2(r,v)} = \epsilon(r)$$

for $1 \leq r \leq n$, $\theta_1(r,v), \theta_2(r,v) \in \mathbb{Z}$, and $\epsilon(r) = \pm 1$. It is known that there is redundancy in the above equations so that exactly $n - k$ of them are independent [9]. We denote the solution set of the above equations in $(\mathbb{C} \backslash \{0,1\})^n$ by $H(M)$ and the point corresponding to the complete structure by $z^0 \in H(M)$.

Let $T_i$ be a torus cross-section of the $i^{\text{th}}$-cusp and $l_i, m_i$ be the chosen longitude-meridian pair of $T_i$ $(1 \leq i \leq k)$. For each $z \in H(M)$, by giving similarity structures on the tori $T_i$, the dilation components of the holonomies (of the similarity structures) of $l_i$ and $m_i$ are represented in the following forms:

$$\delta(z)(l_i) = \pm \prod_{v=1}^{n} z_v^{\lambda_1(i,v)} \cdot (1-z_v)^{\lambda_2(i,v)}$$
(2.2)
$$\delta(z)(m_i) = \pm \prod_{v=1}^{n} z_v^{\mu_1(i,v)} \cdot (1-z_v)^{\mu_2(i,v)}.$$

Then $\delta(z)(l_i)$ and $\delta(z)(m_i)$ behave very nicely near $z^0$ [9].

**Theorem 2.1.** $\delta(z)(l_i) = 1$ and $\delta(z)(m_i) = 1$ are equivalent in a small neighborhood of $z^0$.

**Theorem 2.2.** $z^0$ is a smooth point of $H(M)$ and the unique point near $z^0$ with all $\delta(z)(l_i) = 1$ $(1 \leq i \leq k)$.

By taking logarithms locally near the point $z^0$, equation (2.1) can be re-written as follows:

(2.3) $$\sum_{v=1}^{n} \Big(\theta_1(r,v) \cdot \log(z_v) + \theta_2 \cdot \log(1-z^v)\Big) = c(r) \quad \text{for} \quad r = 1, \ldots, n-k$$

where $c(r)$ are some suitable constants. In the same way, if we let

(2.4) $$u_i(z) = \log\big(\delta(z)(l_i)\big) \quad i = 1, \ldots, k$$
(2.5) $$v_i(z) = \log\big(\delta(z)(m_i)\big) \quad i = 1, \ldots, k$$

in a small neighborhood of $z^0$, then $v_1, \ldots, v_k$ can be parametrized holomorphically in terms of $u_1, \ldots, u_k$ as below [9]:

**Theorem 2.3.** In a neighborhood of the origin in $\mathbb{C}^k$ (with coordinates $u_1, \ldots, u_n$), the following holds for each $i$ $(1 \leq i \leq k)$

(1) $v_i = u_i \cdot \tau_i(u_1, \ldots, u_k)$ where $\tau_i(u_1, \ldots, u_k)$ is an even function of its arguments with $\tau_i(0, \ldots, 0) = \tau_i$ (the cusp shape of $T_i$ with respect to $l_i, m_i$).
(2) There is an analytic function $\Phi(u_1, \ldots, u_k)$ such that $\partial \Phi / \partial u_i = 2v_i$ and $\Phi(0, \ldots, 0) = 0$.
(3) $\Phi(u_1, \ldots, u_k)$ is even in each argument and it has Taylor expansion of the form:

$$\Phi(u_1, \ldots, u_k) = \tau_1 u_1^2 + \cdots + \tau_k u_k^2 + \text{(higher order)}.$$

We call $\Phi(u_1, \ldots, u_k)$ the *potential function* with respect to $u_i, v_i$ $(1 \leq i \leq k)$ and use $\text{Def}(M)$ to denote a small neighborhood of $z^0$ of the manifold defined in (2.3).

Let

$T_p^* \text{Def}(M)$ : The space of holomorphic differentials at $p \in \text{Def}(M)$

(i.e. The cotangent space of type $(1,0)$)

and

$du_i|_p =$ The holomorphic differential induced by $u_i(z)$ at $p \in \text{Def}(M)$

$dv_i|_p =$ The holomorphic differential induced by $v_i(z)$ at $p \in \text{Def}(M)$

where $1 \leq i \leq k$. Then the above theorems imply the following corollary which plays a key role in the proofs of later theorems.

**Corollary 2.4.** *(1) $\{du_1|_p, \ldots, du_k|_p\}$ is a basis of $T_p^*\text{Def}(M)$.*
*(2) $\{dv_1|_p, \ldots, dv_k|_p\}$ is a basis of $T_p^*\text{Def}(M)$.*
*(3) $du_i|_{z^0} = \tau_i dv_i|_{z^0}$ in $T_{z^0}^*\text{Def}(M)$ for $1 \leq i \leq k$.*

**2.2. Holonomy variety (Deformation variety)** There are several ways to define the deformation variety but here we choose the one which is called the holonomy variety, a natural extension of the gluing variety defined in the previous subsection.

Consider the map

$$\xi : z \longrightarrow (\delta(z)(m_1), \ldots, \delta(z)(m_k), \delta(z)(l_1), \ldots, \delta(z)(l_k)) \in \mathbb{C}^{2k},$$

then the *holonomy variety* of $M$ is the Zariski closure $\overline{\xi(H(M))}$ of the image of the above map. In general, the point $(1, \ldots, 1)$ which is the image of the complete structure is a singular point, but there exists a local branch which is isomorphic to $\text{Def}(M)$. From now on, when we say the holonomy variety of a hyperbolic 3-manifold $M$, we indicate the whole variety $\overline{\xi(H(M))}$. But specifically when we mention the irreducible holonomy variety, it only means the irreducible component of it containing the local branch corresponding to $\text{Def}(M)$.

**Remark.** It is a standard fact from algebraic geometry that if a variety is defined over rational numbers, then the Zariski closure of the image of it under a rational map is also defined over rational numbers (thus defining equations of the holonomy variety consists of rational polynomials). Also throughout the paper, the irreducibility means the one over the algebraic closures $\overline{\mathbb{Q}}$ or $\mathbb{C}$.

**2.3. Dehn Surgery** Hyperbolic Dehn surgery (Dehn filling) can be defined in a few slightly different ways. In this paper, we adopt the definition that, after attaching a new torus, the core of the torus is always isotopic to a geodesic of the Dehn filled manifold. This definition will allow us to avoid some redundancy in manifolds and thus simplify the proofs of the main theorems.

Let $M_{p_1/q_1,\ldots,p_k/q_k}$ be the Dehn filled manifold of $M$ with surgery coefficient $(p_1/q_1, \ldots, p_k/q_k)$. By the Seifert-Van Kampen theorem, the fundamental group of $M_{p_1/q_1,\ldots,p_k/q_k}$ is obtained by adding relations

$$m_1^{p_1} l_1^{q_1} = 1, \quad \ldots \quad , m_k^{p_k} l_k^{q_k} = 1$$

to the fundamental group of $M$. Hence, on the holonomy variety of $M$, the hyperbolic structure of $M_{p_1/q_1,\ldots,p_k/q_k}$ is identified with a point satisfying additional equations corresponding to the above relations. More precisely, if the holonomy variety of $M$ is given as

(2.6) $$f_i(M_1, \ldots, M_k, L_1, \ldots, L_k) = 0 \quad (1 \leq i \leq s),$$

then a holonomy representation of $M$ which gives rise to an incomplete structure inducing the Dehn filled manifold $M_{p_1/q_1,\ldots,p_k/q_k}$ is a point satisfying the following equations:

(2.7) $$M_1^{p_1} L_1^{q_1} = 1, \quad \ldots \quad , M_k^{p_k} L_k^{q_k} = 1.$$

We call (2.7) the *Dehn surgery equations* with coefficient $(p_1/q_1, \ldots, p_k/q_k)$ and the two points inducing the hyperbolic structure on $M_{p_1/q_1,\ldots,p_k/q_k}$ the *Dehn filling points* corresponding to $M_{p_1/q_1,\ldots,p_k/q_k}$.

If $m_i^{s_i} l_i^{r_i}$ represents a core curve of the Dehn filled manifold $M_{p_1/q_1,\ldots,p_k/q_k}$ (so that $p_i r_i - q_i s_i = 1$), then the eigenvalue of $m_i^{s_i} l_i^{r_i}$ is of the form $(t_i^{\frac{1}{2}}, t_i^{-\frac{1}{2}})$ or $(-t_i^{\frac{1}{2}}, -t_i^{-\frac{1}{2}})$ where $t_i^{-q_i}$ and $t_i^{p_i}$ are the holonomies of $m_i$ and $l_i$ respectively (i.e. $M_i = t_i^{-q_i}, L_i = t_i^{p_i}$). We let $\epsilon_i(t_1^{\frac{1}{2}} + t_1^{-\frac{1}{2}})$ (where $\epsilon_i = 1$ or $-1$) be the trace value of $m_i^{s_i} l_i^{r_i}$ for each $i$ and name

$$\left(\epsilon_1(t_1^{\frac{1}{2}} + t_1^{-\frac{1}{2}}), \ldots, \epsilon_k(t_k^{\frac{1}{2}} + t_k^{-\frac{1}{2}})\right)$$

the *core trace value* of the Dehn filling coefficient $(p_1/q_1, \ldots, p_k/q_k)$. (Note that $|t_i| \neq 1$ for each $i$ since it's an eigenvalue of a hyperbolic element.)

**Remark.** We could use $(\epsilon_1 t_1, \ldots, \epsilon_k t_k)$ or $(\epsilon_1 t_1^{\frac{1}{2}}, \ldots, \epsilon_k t_k^{\frac{1}{2}})$ instead of $\left(\epsilon_1(t_1^{\frac{1}{2}} + t_1^{-\frac{1}{2}}), \ldots, \epsilon_k(t_k^{\frac{1}{2}} + t_k^{-\frac{1}{2}})\right)$ as it doesn't make any essential difference. But, to the author's perspective, the latter one is more natural and easier to define. It even consists of elements of the trace field of $M$, and so it is more convenient to handle in the proofs.

In the above definition of hyperbolic Dehn surgery, there may exist some Dehn filling points which are not contained in the irreducible holonomy variety. To avoid this issue, we now define a somewhat stronger version of hyperbolic Dehn surgery. More precisely, we say $M_{p_1/q_1,\ldots,p_k/q_k}$ is obtained by strong hyperbolic Dehn filling if its hyperbolic structure can be deformed to the complete structure on $M$ through a family of cone manifolds and all the corresponding points on the representation variety are smooth. Then with this new definition, we can ignore Dehn filling points outside of the irreducible holonomy variety. From now on, when we mention hyperbolic Dehn filling, we always mean this stronger version.

The first theorem below is that of Jorgensen-Thurston which greatly simplifies the structure of hyperbolic 3-manifolds of bounded volume, and the second one is a part of Thurston's hyperbolic Dehn surgery theory (the first one was originally formulated under the previous definition of hyperbolic

Dehn surgery, but it's not hard to see that the theorem is also true under the stronger definition.) [2]:

**Theorem 2.5.** *For any $D > 0$, there exists a finite set of non-compact manifolds $M_1, ..., M_k$ such that all closed hyperbolic 3-manifolds of volume less than or equal to $D$ are obtained by hyperbolic Dehn surgery on $M_i$ for some $i$.*

**Theorem 2.6.** *Using the same notation as above, for each $i$, the value $t_i$ converges to 1 as $max(|p_i|, |q_i|)$ goes to $\infty$.*

## 3. Background II (Number Theory)

3.1. **Mahler measure and length of a polynomial** The Mahler measure $\mathcal{M}(f)$ and length $\mathcal{L}(f)$ of an integer polynomial

$$f(X) = a_n X^n + \cdots + a_1 X + a_0 = a_n(X - \alpha_1) \cdots (X - \alpha_n)$$

are defined by

$$\mathcal{M}(f) = |a_n| \prod_{i=1}^{n} \max(|\alpha_i|, 1),$$

$$\mathcal{L}(f) = |a_0| + \cdots + |a_n|$$

respectively. Then the following properties are standard [6]:

(1) $\mathcal{M}(f_1 f_2) = \mathcal{M}(f_1)\mathcal{M}(f_2)$
(2) $\mathcal{M}(f) \leq \mathcal{L}(f)$

where $f_1$ and $f_2$ are two integer polynomials.

3.2. **Height** The height $H(\alpha)$ of an algebraic number $\alpha$ is defined as follows:

**Definition 3.1.** *Let $K$ be an any number field containing $\alpha$, $M_K$ be the set of places of $K$, and $K_v, \mathbb{Q}_v$ be the completions at $v \in M_K$. Then*

$$H(\alpha) = \prod_{v \in M_K} \max\{1, |\alpha|_v\}^{[K_v:\mathbb{Q}_v]/[K:\mathbb{Q}]}$$

Note that the above definition doesn't depend on the choice $K$. That is, for any number field $K$ containing $\alpha$, it gives us the same value. The following properties can be easily deduced from the definition [4].

**Theorem 3.2.** *(1) There are only finitely many algebraic numbers of uniformly bounded height and degree.*
*(2) $H(\alpha) = H(1/\alpha)$ for $\alpha \in \overline{\mathbb{Q}}$.*
*(3) $H(\alpha_1 + \cdots + \alpha_r) \leq rH(\alpha_1) \cdots H(\alpha_r)$ for $\alpha_1, ..., \alpha_r \in \overline{\mathbb{Q}}$.*
*(4) $H(\alpha_1 \cdots \alpha_r) \leq H(\alpha_1) \cdots H(\alpha_r)$ for $\alpha_1, ..., \alpha_r \in \overline{\mathbb{Q}}$.*
*(5) $H(\alpha)^{\deg f} = \mathcal{M}(f)$ where $\alpha \in \overline{\mathbb{Q}}$ and $f$ is the minimal polynomial of $\alpha$.*

If $\alpha = (\alpha_1, ..., \alpha_n) \in \overline{\mathbb{Q}}^n$ is an $n$-tuple of algebraic numbers, the definition can be generalized as follows:

**Definition 3.3.** *Let $K$ be an any number field containing $\alpha_1, ..., \alpha_n$, $M_K$ be the set of places of $K$, and $K_v, \mathbb{Q}_v$ be the completions at $v$. Then*

$$H(\alpha) = \prod_{v \in M_K} \max\{1, |\alpha_1|_v, ..., |\alpha_n|_v\}^{[K_v:\mathbb{Q}_v]/[K:\mathbb{Q}]}$$

Similar to Theorem 3.2, the following inequalities holds:

(3.1) $\quad \max\{H(\alpha_1), \ldots, H(\alpha_n)\} \leq H(\alpha) \leq H(\alpha_1) \cdots H(\alpha_n).$

Now, using the core trace value, we define the *height of the Dehn filling coefficient* $(p_1/q_1, ..., p_k/q_k)$ by (with the same notations in Section 2.3)

$$H\Big(\big(\epsilon_1(t_1^{\frac{1}{2}} + t_1^{-\frac{1}{2}}), ..., \epsilon_k(t_k^{\frac{1}{2}} + t_k^{-\frac{1}{2}})\big)\Big).$$

We next show that **Question 3** is stronger than **Question 2** by applying Theorem 3.2 (1).

**Theorem 3.4.** *If the answer to **Question 3** is yes, then so is the answer to **Question 2**.*

*Proof.* Suppose that the answer to **Question 3** is yes and $\big(t := \epsilon_1(t_1^{\frac{1}{2}} + t_1^{-\frac{1}{2}}), ..., \epsilon_k(t_k^{\frac{1}{2}} + t_k^{-\frac{1}{2}})\big)$ is the core trace value of an arbitrary Dehn filling coefficient inducing a Dehn filled manifold $M_{dehn}$ of $M$. Clearly $\mathbb{Q}\big(\epsilon_1(t_1^{\frac{1}{2}} + t_1^{-\frac{1}{2}}), ..., \epsilon_k(t_k^{\frac{1}{2}} + t_k^{-\frac{1}{2}})\big)$ (say $\mathbb{Q}(t)$) is a subfield of the trace field of $M_{dehn}$. Since the height of $t$ is bounded by the universal constant, if the degree of $\mathbb{Q}(t)$ is bounded, then there are only finitely many choices for the core trace value $t$ by Theorem 3.2 (1) and (3.1). Furthermore, for the given $t$, there are also only finitely many Dehn surgery coefficients having $t$ as the core trace value because of Theorem 2.2. This completes the proof. $\square$

As observed in Section 2.3, a Dehn filling point inducing the manifold $M_{p_1/q_1,...,p_k/q_k}$ is of the following form:

(3.2) $\quad (M_1, \ldots, M_k, L_1, \ldots, L_k) = (t_1^{-q_1}, \ldots, t_k^{-q_k}, t_1^{p_1}, \ldots, t_k^{p_k}).$

If the height of (3.2) is bounded, then the height of each $t_i$ and the core trace value $\big(\epsilon_1(t_1^{\frac{1}{2}} + t_1^{-\frac{1}{2}}), \ldots, \epsilon_k(t_k^{\frac{1}{2}} + t_k^{-\frac{1}{2}})\big)$ are also bounded by (3.1) and Theorem 3.2. Hence, to prove the uniform boundedness of the heights of the core trace values of Dehn fillings, it is enough to prove the uniform boundedness of the heights of their corresponding Dehn filling points.

**3.3. Anomalous Subvarieties** In this section, we identify $G_m^n$ with the non-vanishing of the coordinates $x_1, \ldots, x_n$ in the affine $n$-space $\overline{\mathbb{Q}}^n$ or $\mathbb{C}^n$ (i.e. $G_m^n = (\overline{\mathbb{Q}}^*)^n$ or $(\mathbb{C}^*)^n$). An algebraic subgroup $H_\Lambda$ of $G_m^n$ is defined as the set of solutions satisfying equations $x_1^{a_1} \cdots x_n^{a_n} = 1$ where the vector $(a_1, \ldots, a_n)$ runs through a lattice $\Lambda \subset \mathbb{Z}^n$. If $\Lambda$ is primitive, then we call $H_\Lambda$ an irreducible algebraic subgroup or algebraic torus. By a coset $K$, we mean a translate $gH$ of some algebraic subgroup $H$ by some $g \in G_m^n$. To simplify notation, for $\mathbf{i} = (i_1, \ldots, i_n) \in \mathbb{Z}^n$, we abbreviate $x_1^{i_1} \cdots x_n^{i_n}$ as $\mathbf{x}^\mathbf{i}$. Let $\mathbf{e}_1 = (1, 0, \ldots, 0)^t, \ldots, \mathbf{e}_n = (0, 0, \ldots, 1)^t$ be column vectors, which we identify with the usual basis of $\mathbb{Z}^n$, and $A$ be an $n \times n$ matrix with columns $A\mathbf{e}_i = (a_{1i}, \ldots, a_{ni})^t \in \mathbb{Z}^n$ for $i = 1, \ldots, n$. Then the map $\varphi_A : G_m^n \longrightarrow G_m^n$ defined by
$$\varphi_A(\mathbf{x}) := (\mathbf{x}^{A\mathbf{e}_1}, \ldots, \mathbf{x}^{A\mathbf{e}_n})$$
is called a monoidal transformation. This is a typical homomorphism of $G_m^n$ and will be repeatedly used throughout the paper. For more properties of algebraic subgroups and $G_m^n$, see [4].

The following theorem is Proposition 3.2.7 in [4]. We include a proof here because the idea behind it will be applied several times later :

**Theorem 3.5.** *Let $H \subset G_m^n$ be an algebraic subgroup of rank $n - r$. Then there exists a monoidal transformation $\phi$ such that $\phi(H)$ is equal to $F \times G_m^r$ ($\subset G_m^{n-r} \times G_m^r = G_m^n$) where $F$ is a finite algebraic subgroup.*

*Proof.* Let $\Lambda$ be a subgroup of $\mathbb{Z}^n$ such that $H_\Lambda = H$. By the theorem of elementary divisors, there is a basis $\mathbf{b_1}, \ldots, \mathbf{b_n}$ of $\mathbb{Z}^n$ and elements $\lambda_1, \ldots, \lambda_{n-r} \in \mathbb{Z} \backslash \{0\}$ such that $\lambda_1 \mathbf{b_1}, \ldots, \lambda_{n-r} \mathbf{b_{n-r}}$ is a basis of $\Lambda$. Using a monoidal transformation to change coefficients, we may assume that $\mathbf{b_1}, \ldots, \mathbf{b_n}$ is the standard basis. Then $H$ is isomorphic to $F \times G_m^r$ with
$$F = \{\mathbf{x} \in G_m^{n-r} \mid x_1^{\lambda_1} = 1, \ldots, x_{n-r}^{\lambda_{n-r}} = 1\}.$$
$\square$

**Definition 3.6.** *An irreducible subvariety $Y$ of $X$ is anomalous (or better, $X$-anomalous) if it has positive dimension and lies in a coset $K$ in $G_m^n$ satisfying*
$$\dim K \leq n - \dim X + \dim Y - 1.$$

The quantity $\dim X + \dim K - n$ is what one would expect for the dimension of $X \cap K$ when $X$ and $K$ were in general position. Thus we can understand anomalous subvarieties of $X$ as the ones that are unnaturally large intersections with cosets of algebraic subgroups of $G_m^n$ (see [13] for more discussions about this).

**Definition 3.7.** *The deprived set $X^{oa}$ is what remains of $X$ after removing all anomalous subvarieties.*

**Definition 3.8.** *An anomalous subvariety of $X$ is maximal if it is not contained in a strictly larger anomalous subvariety of $X$.*

The following theorem tells us the structure of anomalous subvarieties (Theorem 1 of [3]).

**Theorem 3.9.** *Let $X$ be an irreducible variety in $G_m^n$ of positive dimension defined over $\overline{\mathbb{Q}}$.*
*(a) For any torus $H$ with*

(3.3) $$1 \leq h = n - \dim H \leq \dim X$$

*the union $Z_H$ of all subvarieties $Y$ of $X$ contained in any coset $K$ of $H$ with*

(3.4) $$\dim Y = \dim X - h + 1$$

*is a closed subset of $X$, and the product $HZ_H$ is not Zariski dense in $G_m^n$.*
*(b) There is a finite collection $\Psi = \Psi_X$ of such tori $H$ such that every maximal anomalous subvariety $Y$ of $X$ is a component of $X \cap gH$ for some $H$ in $\Psi$ satisfying (3.3) and (3.4) and some $g$ in $Z_H$. Moreover $X^{oa}$ is obtained from $X$ by removing the $Z_H$ of all $H$ in $\Psi$, and thus it is open in $X$ with respect to the Zariski topology.*

Now we recall the bounded height conjecture which we stated in Section 1.

**Theorem 1.2** (Bounded Height Conjecture=Habegger's theorem) *Let $X \subset G_m^n$ be an irreducible variety over $\overline{\mathbb{Q}}$. The height is bounded in the intersection of $X^{oa}$ with the union of algebraic subgroups of dimension $\leq n - \dim X$.*

We next explain how the above theorem fits into the setting of our problem. In our case, $X$, the holonomy variety of a $k$-cusped hyperbolic 3-manifold is a $k$-dimensional variety in the $2k$-dimensional ambient space and Dehn surgery equations define $k$-dimensional algebraic subgroups. So they exactly satisfy the dimension condition cited above. Consequently, as we explained in Section 1, Theorem 1.2 tells us the uniformly boundedness of the heights of the Dehn filling points not on anomalous subvarieties. Hence, to prove the uniformly boundedness of all Dehn filling points, it is enough to analyze the structures of anomalous subvarieties of $X$ and the heights of Dehn filling points on them. Of course, in the worst case, it is possible that $X$ is a maximal anomalous variety of itself and so the Bounded Height Conjecture tells us nothing. But for the holonomy variety of a hyperbolic 3-manifold, we can show that that is not the case.

**Theorem 3.10.** *If $X$ is the irreducible holonomy variety of a $k$-cusped hyperbolic manifold $M$, then $X$ itself is not an anomalous variety.*

*Proof.* If $X$ is anomalous, then $X$ is contained in an algebraic subgroup defined by an equation of the form

$$M_1^{a_1} \cdots M_k^{a_k} L_1^{b_1} \cdots L_k^{b_k} = 1$$

where not all $a_i, b_i$ are zero. Translating this information into $\text{Def}(M)$, it implies the differential of

(3.5) $$a_1 u_1(z) + \cdots + a_k u_k(z) + b_1 v_1(z) + \cdots + b_k v_k(z)$$

at $z^0$, which is

(3.6) $$a_1 du_1|_{z_0} + \cdots + a_k du_k|_{z_0} + b_1 dv_1|_{z_0} + \cdots + b_k dv_k|_{z_0},$$

is zero in $T^*_{z_0}\text{Def}(M)$. By Corollary 2.4, (3.6) is equal to

(3.7) $$(a_1 + b_1\tau_1)du_1|_{z_0} + \cdots + (a_k + b_k\tau_k)du_k|_{z_0}.$$

But (3.7) is zero in $T^*_{z_0}\text{Def}(M)$ iff all the coefficients $a_j, b_j$ are zero since $\tau_j \notin \mathbb{R}$. This contradicts the original assumption on $a_j$ and $b_j$. □

Before going on to the next section, let us briefly go through the 1-cusped case since it provides the basic ideas for the other cases. We prove it using two different methods (i.e. Habegger's theorem and Hodgson's method) as both naturally extend to the higher cusped cases. Originally, Hodgson didn't use the notion of height in his proof, but the key ideas are essentially the same.

**Theorem 3.11.** *For a 1-cusped manifold $M$, there exists a constant $D > 0$ such that the height of any Dehn filled manifold $M_{p/q}$ of $M$ is bounded by $D$.*

*1st Proof.* Let $X$ be the irreducible holonomy variety of $M$. By Definition 4.3, the only possible anomalous subvariety of $X$ is $X$ itself. But this is impossible by Theorem 3.10. So $X^{oa} = X$ and, applying Habegger's theorem, we get the desired result. □

*2nd Proof.* Let $f(M, L) = 0$ be the defining equation of the holonomy variety with integer coefficients. If a Dehn filling equation is given by $M^p L^q = 1$, a corresponding Dehn filling point is of the form $M = t^{-q}, L = t^p$. By multiplying by a power of $t$ if needed, we may assume $f(t^{-q}, t^p)$ is an integer polynomial. Then the following inequalities hold by (1), (2) in Section 3.1 and Theorem 3.2:

$$H(t) \leq \mathcal{M}(f(t^{-q}, t^p)) \leq \mathcal{L}(f(t^{-q}, t^p)) \leq \mathcal{L}(f(M, L))$$

where $\mathcal{L}(f(M, L))$ is the sum of the absolute values of all the coefficients of $f(M, L)$. This implies

$$H(t + 1/t) \leq 2H(t)H(1/t) = 2H(t)^2 \leq 2\mathcal{L}\big(f(M, L)\big)^2.$$

Hence all the height of any Dehn filling point is uniformly bounded by $2\mathcal{L}\big(f(M, L)\big)^2$. □

## 4. 2-cusped case

By Theorem 3.10, since the irreducible holonomy variety $X$ is not itself anomalous, we have the following dichotomy for the 2-cusped case:

**Type (I)** $X$ has only a finite number of maximal anomalous subvarieties.

**Type (II)** $X$ has an infinite number of maximal anomalous subvarieties. More specifically, there exists an algebraic subgroup $H$ so that $X$ is foliated by subvarieties contained in $\bigcup_{g \in Z_H} gH \cap X$.

Surely, in the case of **Type (II)**, we cannot apply Habegger's theorem because $X^{oa} = \emptyset$. But soon we will see that this is closely related to a certain geometric phenomenon, namely strong geometric isolation, mentioned in Section 1. In this case, as explained in the same section, the uniformly boundedness of the heights of Dehn filling points follows by extending Hodgson's method.

4.1. **Strong Geometric Isolation** Strong geometric isolation was first introduced by W. Neumann and A. Reid in [8]. Geometrically, this simply means that one subset of cusps moves independently without affecting the rest. Using Theorem 4.3 in [8], we give one of the equivalent forms of the definition as follows:

**Definition 4.1.** *Let $M$ be a $k$-cusped hyperbolic 3-manifold. We say cusps $1, \ldots, l$ are strongly geometrically isolated from cusps $l+1, \ldots, k$ if $v_1, \ldots, v_l$ only depend on $u_1, \ldots, u_l$ and not on $u_{l+1}, \ldots, u_k$.*

When a manifold has this property for each cusp, i.e. each $v_k$ depends only on $u_k$, then the uniformly boundedness holds by Hodgson's method.

**Theorem 4.2.** *If $M$ is a $k$-cusped hyperbolic 3-manifold whose cusps are strongly geometrically isolated from each other, then the height of any Dehn filling point of its irreducible holonomy variety $X$ is uniformly bounded.*

*Proof.* By Definition 4.1, every holonomy $v_i$ is a function of the single variable $u_i$ where $1 \leq i \leq k$. Using the same notation as in Theorem 2.3, we have $v_i = u_i \tau_i(u_i)$. For each $i$ ($1 \leq i \leq k$), consider the following projection

$$\xi_i : X \longrightarrow \mathbb{C}^2$$
$$(M_1, \ldots, L_k) \longmapsto (M_i, L_i).$$

Then $\overline{\xi_i(X)}$ is an algebraic curve which contains a local branch isomorphic to $v_i = u_i \tau_i(u_i)$. Let $f_i(M_i, L_i) = 0$ be a defining polynomial of $\overline{\xi_i(X)}$ having integer coefficients. Then the variety defined by $f_i(M_i, L_i) = 0$ ($1 \leq i \leq k$) is a $k$-dimensional variety containing $X$.

By the second proof of Theorem 3.11, the height of any Dehn filling point of $X$ is bounded by the maximum of $\{2\mathcal{L}(f_i(M_i, L_i))^2 : 1 \leq i \leq k\}$. This completes the proof. $\square$

4.2. **Rationally Independent Cusps with Infinitely Many Anomalous Subvarieties** To prove Theorem 1.4, we first prove the following lemma, which is of independent interest.

**Lemma 4.3.** *Let $M$ be a 2-cusped manifold with rationally independent cusp shapes. Then the only maximal anomalous varieties of its irreducible* holonomy variety $X$ containing $(1,1,1,1)$ are cut out by $M_1 = 1, L_1 = 1$ and $M_2 = 1, L_2 = 1$.

*Proof.* By Theorem 3.9, any maximal anomalous subvariety $Y$ of $X$ containing $(1,1,1,1)$ is of the following form: there exists a 2-dimensional algebraic torus $H$ such that $Y \subset H \cap X$ and $\dim Y = 1$ (this is the only case satisfying the dimension conditions (3.3) and (3.4) in Theorem 3.9).

Let $H$ be given by

$$M_1^{a_1} L_1^{b_1} M_2^{c_1} L_2^{d_1} = 1$$
$$M_1^{a_2} L_1^{b_2} M_2^{c_2} L_2^{d_2} = 1.$$

Moving to $\mathrm{Def}(M)$, if we set

$$h_i = a_i u_1(z) + b_i v_1(z) + c_i u_2(z) + d_i v_2(z)$$
$$dh_i|_{z^0} = \text{The holomorphic differential of } h_i \text{ at } z^0 \ (i = 1, 2),$$

then

(4.1) $\quad dh_i|_{z^0} = (a_i + b_i \tau_1) du_1|_{z^0} + (c_i + d_i \tau_2) du_2|_{z^0} \ (i = 1, 2)$

in $T_{z^0}^* \mathrm{Def}(M)$ by Corollary 2.4 (where $\tau_1, \tau_2$ are the cusp shapes as usual). Since $\dim H \cap X = 1$, the dimension of the space $\langle dh_1|_{z^0}, dh_2|_{z^0} \rangle$ is also equal to 1 in $T_{z^0}^* \mathrm{Def}(M)$. (If $z^0$ is a singular point of $H \cap X$, then $\dim \langle dh_1|_{z^0}, dh_2|_{z^0} \rangle$ could be 0. But, in this case, we get $a_i = b_i = c_i = d_i = 0$ for $i = 1, 2$, which contradicts our original assumption that $H$ is a 2-dimensional algebraic subgroup.)

By (4.1), it can be shown that $\dim \langle dh_1|_{z^0}, dh_2|_{z^0} \rangle = 1$ in $T_{z^0}^* \mathrm{Def}(M)$ iff

$$(a_1 + b_1 \tau_1)(c_2 + d_2 \tau_2) = (c_1 + d_1 \tau_1)(a_2 + b_2 \tau_2),$$

and as $1, \tau_1, \tau_2, \tau_1 \tau_2$ are linearly independent over $\mathbb{Q}$, this is equivalent to

(4.2) $\quad\quad\quad\quad\quad\quad a_1 c_2 - c_1 a_2 = 0$

(4.3) $\quad\quad\quad\quad\quad\quad b_1 c_2 - c_1 b_2 = 0$

(4.4) $\quad\quad\quad\quad\quad\quad a_1 d_2 - d_1 a_2 = 0$

(4.5) $\quad\quad\quad\quad\quad\quad b_1 d_2 - d_1 b_2 = 0.$

**Claim 4.4.** *The equations (4.2)-(4.5) induce either $a_1 = a_2 = b_1 = b_2 = 0$ (with $c_1 d_2 - c_2 d_1 \neq 0$) or $c_1 = c_2 = d_1 = d_2 = 0$ (with $a_1 b_2 - a_2 b_1 \neq 0$).*

*Proof.* If none of $a_i, b_i, c_i, d_i$ $(i = 1, 2)$ are zero, then (4.2)-(4.5) imply the two nonzero vectors $(a_1, b_1, c_1, d_1)$ and $(a_2, b_2, c_2, d_2)$ are linearly dependent over $\mathbb{Q}$. But this is impossible because $H$ is a 2-dimensional algebraic subgroup. Without loss of generality, let's assume $a_1 = 0$. Then, by (4.2) and (4.4), we have the following two cases:

**Case 1.** $a_2 = 0$

In this case, the problem is reduced to the following:

(4.6) $$b_1 c_2 - c_1 b_2 = 0,$$
(4.7) $$b_1 d_2 - d_1 b_2 = 0,$$
(4.8) $\quad (b_1, c_1, d_1)$ and $(b_2, c_2, d_2)$ are linearly independent.

Just like above, if none of $b_i, c_i, d_i$ ($i = 1, 2$) are zero, then $(b_1, c_1, d_1)$ and $(b_2, c_2, d_2)$ are linearly dependent over $\mathbb{Q}$ by (4.6) and (4.7), contradicting (4.8). So at least one of $b_i, c_i, d_i$ ($i = 1, 2$) is zero and the situation is divided into the following two subcases.

(1) $b_1 = 0$ or $b_2 = 0$

By symmetry, it is enough to consider the case $b_1 = 0$. If $b_1 = 0$, then $b_2 = 0$ or $c_1 = 0$ (from (4.6)) and $b_2 = 0$ or $d_1 = 0$ (from (4.7)). If $b_2 = 0$, then we get the desired result (i.e. $a_1 = a_2 = b_1 = b_2 = 0$). Otherwise, if $c_1 = d_1 = 0$, this contradicts the fact that $(a_1, b_1, c_1, d_1)$ is a nonzero vector.

(2) $c_1 = 0$ or $c_2 = 0$ or $d_1 = 0$ or $d_2 = 0$ (with $b_1, b_2 \neq 0$)

Here, also by symmetry, it is enough to prove the first case $c_1 = 0$. If $b_1, b_2 \neq 0$ and $c_1 = 0$, then $c_2 = 0$ by (4.6) and the problem is further simplified to the following:

$$b_1 d_2 - d_1 b_2 = 0,$$

$(b_1, d_1)$ and $(b_2, d_2)$ are linearly independent.

However this doesn't hold regardless of the values of $d_1$ and $d_2$.

**Case 2.** $a_2 \neq 0$ and so $c_1 = d_1 = 0$.

Since $(a_1, b_1, c_1, d_1)$ is a nonzero vector, $b_1$ is nonzero and $c_2 = d_2 = 0$ by (4.3) and (4.5). As a result, we get $c_1 = c_2 = d_1 = d_2 = 0$, which is the second desired result of the statement.

So Claim 4.4 holds. $\square$

We now use Claim 4.4 to complete the proof of Lemma 4.3. Without loss of generality, let's assume $a_1 = a_2 = b_1 = b_2 = 0$ and $H$ is defined by

(4.9)
$$M_2^{c_1} L_2^{d_1} = 1$$
$$M_2^{c_2} L_2^{d_2} = 1.$$

Since $c_1 d_2 - c_2 d_1 \neq 0$, both $M_2$ and $L_2$ are roots of unity. As $H$ is an algebraic torus containing $(1, 1, 1, 1)$, equation (4.9) is equal to $M_2 = 1, L_2 = 1$.

In the same way, one gets $M_1 = 1, L_1 = 1$ from the other assumption $c_1 = c_2 = d_1 = d_2 = 0$. This completes the proof of Lemma 4.3. $\square$

Now we prove the following theorem which implies Theorem 1.4.

**Theorem 4.5.** *Let $M$ be a 2-cusped manifold with its irreducible holonomy variety $X$ and suppose that $M_1 = 1, L_1 = 1$ and $M_2 = 1, L_2 = 1$ are the only algebraic subgroups generating maximal anomalous subvarieties of $X$* *containing* $(1, 1, 1, 1)$. *If $X$ has an infinite number of maximal anomalous subvarieties, then the two cusps of $M$ are strongly geometrically isolated.*

*Proof.* Since $X$ has a maximal anomalous subvariety of **Type (III)**, a small neighborhood $N$ of $(1, 1, 1, 1)$ in $X$ is covered by the intersections of itself with cosets of $M_1 = 1, L_1 = 1$ or $M_2 = 1, L_2 = 1$. Without loss of generality, let's assume the first case. Moving to $\text{Def}(M)$, the given information implies, for any $p = (m_1, l_1, m_2, l_2) \in \text{Def}(M)$, the intersection of $u_1 = m_1, v_1 = l_1$ with $\text{Def}(M)$ is a 1-dimensional submanifold of $\text{Def}(M)$ and so $\dim \langle dv_1|_p, du_1|_p \rangle = 1$ in $T_p^*\text{Def}(M)$. Let $v_1 = g(u_1, u_2)$ and

$$dv_1|_p = g_{u_1}(m_1, m_2)du_1|_p + g_{u_2}(m_1, m_2)du_2|_p.$$

Then $g_{u_2}(m_1, m_2) = 0$ by the aforementioned dimension condition. Since $p$ was chosen arbitrary from $\text{Def}(M)$, one concludes $g_{u_2}(u_1, u_2) = 0$. That is, $v_1$ is a single variable function of $u_1$. If we set $v_2 = h(u_1, u_2)$, then, as $g_{u_2}(u_1, u_2) = h_{u_1}(u_1, u_2)$ (Theorem 2.3), $v_2$ is also a function of single variable $u_2$. This completes the proof. □

**Theorem 1.4** *Let $M$ be a 2-cusped manifold with rationally independent cusp shapes. If the irreducible holonomy variety $X$ of $M$ has an infinite number of maximal anomalous subvarieties, then the two cusps of $M$ are strongly geometrically isolated.*

*Proof.* This immediately follows from Lemma 4.3 and Theorem 4.5. □

## 5. How to Approach to the General Case?

**5.1. Generalized Bounded Height Conjecture** You might wonder why we've been dealing with the holonomy variety instead of the gluing variety even though the holonomy variety is derived from the gluing variety. The reason is that, first, the gluing variety can be a maximal anomalous subvariety of itself. As given in (2.3), the defining equations of the gluing variety are of the following forms:

$$(5.1) \quad \prod_{v=1}^{n} z_v^{\theta_1(r,v)} (1 - z_v)^{\theta_2(r,v)} = \pm 1 \quad (1 \leq r \leq n - k, \ 1 \leq v \leq n).$$

For any fixed $r$, if $\theta_2(r, v) = 0$ for all $v$, then (5.1) contains an equation of an algebraic subgroup so that the gluing variety itself becomes anomalous. Also, in this context, the equations corresponding to Dehn filling equations are represented by equations of the following form (see (2.2))

$$\bigl(\delta(z)(l_i)\bigr)^{p_i} \bigl(\delta(z)(m_i)\bigr)^{q_i} = 1 \quad (1 \leq i \leq k),$$

but these may not define an algebraic subgroup. While there might exist some way to avoid these issues by choosing $z_v$ very delicately, there's no canonical way to do so. Clearly these facts indicate the gluing variety is not a good scheme to work with the Bounded Height Conjecture.

However the holonomy variety itself also has its own flaws. Since the point $(1, \ldots, 1)$ is singular, some pathological things may happen, which prevent us to apply Habegger's theorem. To explain this in more detail, we first assume the holonomy variety has only a finite number of maximal anomalous subvarieties. Let $N_1$ be the branch containing all Dehn filling points (i.e. the one isomorphic to $\mathrm{Def}(M)$) and $N_2$ be another branch intersecting $N_1$ through $(1, \ldots, 1)$ (as $(1, \ldots, 1)$ is singular, it is a reasonable assumption). We further suppose that there exists an algebraic torus $H$ which intersects $N_2$ anomalously but intersects $N_1$ transversally.[1] So far when we worked on these problems in the previous sections, we moved the target from the holonomy variety to $\mathrm{Def}\,(M)$ and used some properties of an anomalous intersection such as the dimension of the cotangent space (for example, see the proof of Lemma 4.3). But in the above setting, if $H$ intersects $N_1$ transversally, it's very hard to characterize the properties of $H \cap N_1$. In particular, to get the uniform boundedness of heights, we need to control the behaviors of the Dehn filling points on $H \cap N_2$ (e.g. whether there are finite or not). However we don't have any information to do this because we know nothing about $N_2$. Also when the holonomy variety has a finite number of maximal anomalous subvarieties, one clever way to attack the problem is applying the Bounded Height Conjecture again to each maximal anomalous subvariety. But, in this case, we will see later that it's impossible to do this because it does not satisfy the required dimension condition of the theorem (see the remark after the proof of Theorem 1.3).

However there's still hope that we can resolve this problem. Intuitively, in the statement of Habegger's theorem, it's very likely that the height may be unbounded only on the branch that intersects with algebraic subgroups anomalously. In other words, although the height is unbounded on $H \cap N_2$ (more precisely, in the intersection of it with the union of algebraic subgroups), we may still expect that it is bounded on $H \cap N_1$. Expanding this idea, we approach the problem in a slightly different way. Thinking of the holonomies of the meridian-longitude pairs $\delta(z)(l_j)$ and $\delta(z)(m_j)$ $(1 \leq j \leq k)$ in (2.2) as new variables, we assume the gluing variety is defined by the following equations in $\mathbb{C}^{2k+n}$.

(5.2)
$$\prod_{v=1}^{n} z_v^{\theta_1(r,v)} \cdot (1-z_v)^{\theta_2(r,v)} = \epsilon(r)$$
$$L_j = \pm \prod_{v=1}^{n} z_v^{\lambda_1(j,v)} \cdot (1-z_v)^{\lambda_2(j,v)}$$
$$M_j = \pm \prod_{v=1}^{n} z_v^{\mu_1(j,v)} \cdot (1-z_v)^{\mu_2(j,v)}$$

---

[1]Here we consider $N_1$, $N_2$, and $H$ as analytic sets

where $1 \leq r \leq n-k$ and $1 \leq j \leq k$. (Similar to $H(M)$ in Section 2.1, we only consider the points satisfying $z_v \neq 0, 1$ for each $v$.) Then the point corresponding to the complete structure is still smooth and the holonomy variety is the Zariski closure of the image of the following natural projection

$$\mathrm{Pr_{2k}} : (L_1, \ldots, M_k, z_1, \ldots, z_n) \longrightarrow (L_1, \ldots, M_k).$$

For convenience, let $X'$ be the gluing variety defined by (5.2) and $X$ ($= \overline{\mathrm{Pr_{2k}}(X')}$) be the holonomy variety. If there exists an anomalous subvariety $Z$ of $X$ produced by the nongeneric intersection between $N_1$ and an algebraic subgroup (where $N_1$ is the same notation in the preceding paragraph), we lift this up, getting the corresponding anomalous subvariety $Z'$ in $X'$. More precisely, if $H \subset (\overline{\mathbb{Q}}^*)^{2k}$ is an algebraic subgroup such that $H \cap X$ contains $Z$, then we pull $H$ back, getting an algebraic subgroup $H' = \mathrm{Pr}_{2k}^{-1}(H) \subset (\overline{\mathbb{Q}}^*)^{2k+n}$ such that $H' \cap X'$ contains $Z'$. Here the point is that both $H$ and $H'$ are defined by exactly the same equations. Moreover, as Dehn filling equations are also defined by the same form (i.e. (2.7)) in the both spaces $(\overline{\mathbb{Q}}^*)^{2k+n}$ and $(\overline{\mathbb{Q}}^*)^{2k}$, we can expect something a similar result by modifying the conditions of the Bounded Height Conjecture slightly. Indeed it turns out that, by following the original proof of Habegger's paper, we can actually prove a generalized version of the theorem that is exactly what we need in our situation, which is summarized as follows:

**Definition 5.1.** *If an algebraic subgroup $H \subset G_m^{n+t}$ contains $\{1\}^n \times G_m^t$, we say that $H$ is defined by restricting the first n-coordinates.*

**Theorem 5.2.** *Let $X \subset G_m^{n+t}$ be an s-dimensional variety ($s \leq n$) and $X^{re-oa}$ be the deprived set after removing all the anomalous subvarieties of $X$ produced by algebraic subgroups (and their cosets) which are defined by restricting the first n-coordinates. Then the height is bounded in the intersection of $X^{re-oa}$ with the union of algebraic subgroups defined by restricting the first n-coordinates and of codimension at least s.*

The proof of the above theorem requires purely number theoretic arguments, so we postpone it until the last section.

5.2. **General case (Theorem 1.3)** Before we prove the main theorem, we first cite two theorems which are necessary to prove it. The first theorem follows from Thurston's hyperbolic Dehn filling theory combined with Mostow's rigidity theorem and the second theorem is Proposition 3.28 in [7].

**Theorem 5.3.** *Let $X$ be the gluing variety of a hyperbolic 3-manifold $M$ and $p$ be a Dehn filling point on $X$. If $K$ is an algebraic subgroup defined by the Dehn filling equations corresponding to $p$, then $p$ is an isolated point in $X \cap K$.*

**Theorem 5.4.** *Let $Z$ be an affine variety and $X, Y$ be subvarieties of $Z$. Let $x \in X \cap Y$ and assume $x$ is smooth on $Z$. If $W$ is any irreducible component* *of $X \cap Y$ containing $x$, then*
$$\dim W \geq \dim X + \dim Y - \dim Z.$$

The following definition confines our targets to the ones having a nice property.

**Definition 5.5.** *We define simple varieties in $G_m^n$ as follows.*

*i) Every 1-dimensional variety is a simple variety.*

*ii) Every algebraic coset is a simple variety.*

*iii) Let $X \subset G_m^n$ be a k-dimensional irreducible variety ($k \geq 2$) and suppose that $X$ is not contained in a proper algebraic subgroup of $G_m^n$. If $X$ has only a finite number of maximal anomalous subvarieties which are all simple, then $X$ is a simple variety.*

*iv) Let $X \subset G_m^n$ be an irreducible variety contained in a proper algebraic subgroup (or coset) of $G_m^n$ and $H$ be an algebraic torus (or coset) of minimal dimension (say $s$) containing $X$. Then we say $X$ is simple if it is simple when regarded as a subvariety in $G_m^s$ ($\cong H$).*

The idea behind this definition is fairly simple. Let $X$ be a simple variety, $Y$ be an anomalous subvariety of $X$ and $H$ be an algebraic torus of minimal dimension such that $Y \subset H \cap X$. Then we get either $Y = H$ or $Y$ has only a finite number of maximal anomalous subvarieties in $H$.

In the following theorem, we assume that the irreducible holonomy variety is simple, but it doesn't necessary mean that the gluing variety is also simple. However, locally near a point corresponding to the complete structure, it contains only a finite number of anomalous subvarieties, which are cut out by algebraic groups containing $\{1\}^{2k} \times G_m^t$, and this is what will matter in the proof.

Now we are ready to prove Theorem 1.3.

**Theorem 1.3** *Suppose that the answer is yes to Question 3 for any s-cusped manifolds where $1 \leq s \leq k-1$. Let $X$ be the irreducible holonomy variety of a k-cusped hyperbolic 3-manifold $M$. If $X$ is simple, then the height of any Dehn filling point on $X$ is uniformly bounded.*

*Proof.* By Thurston's hyperbolic Dehn filling theory and the assumption given in the first sentence, it is enough to prove that the heights of Dehn filling points in an arbitrary small neighborhood of $(1, \ldots, 1)$ is uniformly bounded. (There may exist an infinite number of Dehn filling points outside of this small neighborhood; however all other Dehn filling points are fillings on a finite list of manifolds with $< k$ cusps.)

By shrinking the size of a neighborhood of $(1, \ldots, 1)$ if necessary, we assume that all the maximal anomalous subvarieties of $X$ contain $(1, \ldots, 1)$.

From now on, we work on the gluing variety which we introduced in Section 5.1. By abuse of notation, let's still denote this variety by $X$. Then, near a small neighborhood of the point corresponding to the complete structure, $X$ has a finite number of maximal anomalous subvarieties produced by algebraic subgroups defined by restricting the first $2k$-coordinates.

Suppose $p \in X$ is an arbitrary Dehn filling point and $K$ an algebraic subgroup defined by the Dehn filling equations of $p$. If $p \in X^{re-oa}$, then one gets the desired result by Theorem 5.2. So we assume $p$ is contained in a maximal anomalous subvariety $Y$ of $X$ and $H$ is an algebraic torus containing $\{1\}^{2k} \times (\overline{\mathbb{Q}}^*)^n$ and of the minimal dimension satisfying $Y \subset H \cap X$. Let $\dim Y = l$ and $\dim H = h$.

**Case 1.** If $h = l$, then $Y = H \cap X = H$ (i.e. $H \subset X$). By Theorem 5.3, $p$ is an isolated point in $K \cap X$ and this implies $\dim H \cap K = 0$. So each coordinate of $p$ is a root of unity, but this contradicts the fact that $p$ is a Dehn filling point.

**Case 2.** If $h > l$, then, by Theorem 3.5, there exists a monoidal transformation $\varphi$ such that $\varphi(H) = \{1\}^{2k-h} \times (\overline{\mathbb{Q}}^*)^{h-n} \times (\overline{\mathbb{Q}}^*)^n$. Here, in choosing $\varphi$, we suppose $\varphi$ is the identify map on $\{1\}^{2k} \times (\overline{\mathbb{Q}}^*)^n$ (as $H$ is an algebraic subgroup containing $\{1\}^{2k} \times (\overline{\mathbb{Q}}^*)^n$, this is surely possible). For convenience, we simply say $\varphi(H) = (\overline{\mathbb{Q}}^*)^{h-n} \times (\overline{\mathbb{Q}}^*)^n$ and $\varphi(Y) \subset (\overline{\mathbb{Q}}^*)^{h-n} \times (\overline{\mathbb{Q}}^*)^n$. Because the irreducible holonomy variety is simple and, near the point corresponding to the complete structure, one of its maximal anomalous subvarieties is locally isomorphic to $\varphi(Y)$, $\varphi(Y)$ has only a finite number of maximal anomalous subvarieties produced by algebraic subgroups containing $\{1\}^{h-n} \times (\overline{\mathbb{Q}}^*)^n$.

(1) $\varphi(p) \in \varphi(Y)^{re-oa}$

As $p$ is an isolated point in $K \cap X$ (Theorem 5.3), $p$ is an isolated point in $K \cap Y$ ($\subset K \cap X$) as well. Applying Theorem 5.4 to $X := H \cap K$, $Y := Y$ and $Z := H$, we get $\dim H \cap K \leq h - l$. Since $\varphi(H \cap K)$ is an algebraic subgroup containing $\{1\}^{h-n} \times (\overline{\mathbb{Q}}^*)^n$ in $(\overline{\mathbb{Q}}^*)^{h-n} \times (\overline{\mathbb{Q}}^*)^n$, by Theorem 5.2, it follows that the height of $\varphi(p) \in \varphi(Y)^{re-oa} \cap \varphi(H \cap K)$ is uniformly bounded.

(2) If $\varphi(p) \notin \varphi(Y)^{re-oa}$, we repeat the above process. That is, we find an algebraic subgroup (say $H'$) which produces the maximal anomalous variety of $\varphi(Y)$ containing $\varphi(p)$, then project onto $H'$ and working there by following exactly the same steps. Since the dimension of a new maximal anomalous subvariety decreases whenever we repeat the process, the whole procedures terminate in finitely many steps. So we get the desired result.

□

**Remark.** In the proof of the above theorem, the key point, which enabled us to apply the generalized bounded height conjecture repeatedly, was Theorem 5.3. However, if $X$ is the holonomy variety, then Theorem 5.3 may not be true. For instance, using the same notation as in Section 5.1, if there exists a Dehn filling point $p$ contained in $N_1 \cap N_2$, then $p$ is an isolated point in $N_1 \cap K$ but not necessary in $N_2 \cap K$ (where $K$ is the algebraic subgroup defined by Dehn filling equations corresponding to $p$). As a result, $p$ may not be an isolated point in $X \cap K$ (though this seems to be very unlikely).

## 6. Example

In [1], J. Aaber and N. Dunfield studied the complement of the $(-2, 3, 8)$-pretzel link, which is the sibling of the Whitehead link complement. We denote this hyperbolic manifold by $W$. In their paper, the coefficients of the potential function $\Phi(u_1, u_2)$ of $W$ up to homogeneous degree of 4 were computed, and hence $v_1 = \frac{1}{2}\partial\Phi/\partial u_1$ and $v_2 = \frac{1}{2}\partial\Phi/\partial u_2$ are given as follows:

$$
\begin{aligned}
v_1 &= iu_1 + \left(\frac{-3+i}{48}\right) u_1^3 - \left(\frac{1+i}{16}\right) u_1 u_2^2 + \cdots \\
v_2 &= iu_2 + \left(\frac{-3+i}{48}\right) u_2^3 - \left(\frac{1+i}{16}\right) u_2 u_1^2 + \cdots .
\end{aligned}
\tag{6.1}
$$

Note that since two cuspshapes are the same in this example, we cannot apply Theorem 1.1. Using (6.1), we now prove Theorem 1.5.

**Theorem 1.5** *The irreducible holonomy variety of $W$ is simple (i.e. it has only a finite number of anomalous subvarieties).*

*Proof.* To the contrary, assume the irreducible holonomy variety $X$ of $W$ contains an infinite number of anomalous subvarieties. Let $H$ be a 2-dimensional algebraic subgroup which, along with its cosets, produces infinitely many anomalous subvarieties and

$$
\begin{aligned}
M_1^{a_1} L_1^{b_1} M_2^{c_1} L_2^{d_1} &= 1 \\
M_1^{a_2} L_1^{b_2} M_2^{c_2} L_2^{d_2} &= 1,
\end{aligned}
$$

be the defining equations of $H$. Since $H$ is always contained in an algebraic subgroup defined by equations of the following forms

$$
\begin{aligned}
M_1^{a_1'} L_1^{b_1'} M_2^{c_1'} L_2^{d_1'} &= 1 \\
M_1^{a_2'} L_1^{b_2'} M_2^{c_2'} &= 1,
\end{aligned}
$$

without loss of generality, we may assume $d_2 = 0$. Moving to Def $(W)$, the intersection of any translate of

$$
\begin{aligned}
a_1 u_1 + b_1 v_1 + c_1 u_2 + d_1 v_2 &= 0 \\
a_2 u_1 + b_2 v_1 + c_2 u_2 &= 0.
\end{aligned}
$$

with Def $(W)$ is a 1-dimensional submanifold of Def$(W)$. So the dimension of the space generated by

(6.2)
$$a_1 du_1|_p + b_1 dv_1|_p + c_1 du_2|_p + d_1 dv_2|_p,$$
$$a_2 du_1|_p + b_2 dv_1|_p + c_2 du_2|_p$$

is equal to 1 in $T_p^*\text{Def}(W)$ for any $p \in \text{Def}(W)$. Let $v_1 = g(u_1, u_2)$, $v_2 = h(u_1, u_2)$ and

(6.3)
$$dv_1 = g_{u_1}(u_1, u_2)du_1 + g_{u_2}(u_1, u_2)du_2,$$
$$dv_2 = h_{u_1}(u_1, u_2)du_1 + h_{u_2}(u_1, u_2)du_2.$$

By plugging (6.3) into (6.2), we get

$$\Big(a_1 + b_1 g_{u_1}(m_1, m_2) + d_1 h_{u_1}(m_1, m_2)\Big)du_1|_p + \Big(c_1 + b_1 g_{u_2}(m_1, m_2) + d_1 h_{u_2}(m_1, m_2)\Big)du_2|_p,$$
$$\Big(a_2 + b_2 g_{u_1}(m_1, m_2)\Big)du_1|_p + \Big(c_2 + b_2 g_{u_2}(m_1, m_2)\Big)du_2|_p,$$

which span a 1-dim vector space in $T_p^*\text{Def}(W)$ where $p = (m_1, m_2, g(m_1, m_2), h(m_1, m_2))$. Then, since $p$ is arbitrary, this induces the following equality:

(6.4)
$$\Big(a_1 + b_1 g_{u_1}(u_1, u_2) + d_1 h_{u_1}(u_1, u_2)\Big)\Big(c_2 + b_2 g_{u_2}(u_1, u_2)\Big)$$
$$= \Big(c_1 + b_1 g_{u_2}(u_1, u_2) + d_1 h_{u_2}(u_1, u_2)\Big)\Big(a_2 + b_2 g_{u_1}(u_1, u_2)\Big).$$

Using (6.1), equation (6.4) can be expanded as follows:

(6.5)
$$\left(a_1 + b_1\left(i + \left(\frac{-3+i}{16}\right)u_1^2 - \left(\frac{1+i}{16}\right)u_2^2 + \cdots\right) + d_1\left(-\left(\frac{1+i}{8}\right)u_1 u_2 + \cdots\right)\right)$$
$$\left(c_2 + b_2\left(-\left(\frac{1+i}{8}\right)u_1 u_2 + \cdots\right)\right)$$
$$= \left(c_1 + b_1\left(-\left(\frac{1+i}{8}\right)u_1 u_2 + \cdots\right) + d_1\left(i + \left(\frac{-3+i}{16}\right)u_2^2 - \left(\frac{1+i}{16}\right)u_1^2 + \cdots\right)\right)$$
$$\left(a_2 + b_2\left(i + \left(\frac{-3+i}{16}\right)u_1^2 - \left(\frac{1+i}{16}\right)u_2^2 + \cdots\right)\right).$$

Comparing the coefficients of the constant function, $u_1^2$ and $u_2^2$ in (6.5), we get

(6.6) $$(a_1 + ib_1)c_2 = (c_1 + id_1)(a_2 + ib_2),$$

(6.7) $$-\frac{1+i}{16}d_1(a_2 + ib_2) + \frac{-3+i}{16}b_2(c_1 + id_1) = \frac{-3+i}{16}b_1 c_2,$$

(6.8) $$\frac{-3+i}{16}d_1(a_2 + ib_2) - \frac{1+i}{16}b_2(c_1 + id_1) = -\frac{1+i}{16}b_1 c_2,$$

and hence

(6.9) $$a_1 c_2 = c_1 a_2 - b_2 d_1,$$

(6.10) $$b_1 c_2 = b_2 c_1 + a_2 d_1,$$

(6.11) $$-\frac{1+i}{16} d_1(a_2 + ib_2) = \frac{-3+i}{16}(b_1 c_2 - b_2 c_1 - id_1 b_2),$$

(6.12) $$\frac{-3+i}{16} d_1(a_2 + ib_2) = -\frac{1+i}{16}(b_1 c_2 - b_2 c_1 - ib_2 d_1).$$

Combining (6.10) with (6.11) and (6.12), it follows that

(6.13)
$$-\frac{1+i}{16} d_1(a_2 + ib_2) = \frac{-3+i}{16} d_1(a_2 - ib_2),$$
$$\frac{-3+i}{16} d_1(a_2 + ib_2) = -\frac{1+i}{16} d_1(a_2 - ib_2).$$

Now it is easy to check that (6.13) forces $a_2 = b_2 = 0$ or $d_1 = 0$.

First, if $a_2 = b_2 = 0$, then $c_2 \neq 0$ (otherwise it contradicts the definition of $H$) and $a_1 = b_1 = 0$ from (6.6). As a result, the defining equations of $H$ can be further simplified to $M_2 = 1, L_2 = 1$. But in this case, by Theorem 4.5, $W$ must be a strongly geometrically isolated manifold, which is not true.

Second, if $d_1 = 0$, then the definition of $H$ and (6.6) induces the following fact:

$$c_1 a_2 - a_1 c_2 = 0,$$
$$c_1 b_2 - b_1 c_2 = 0,$$

$(a_1, b_1, c_1)$ and $(a_2, b_2, c_2)$ are linearly independent.

As we observed in Claim 4.4 (Case 1), the only possible case is $c_1 = c_2 = 0$. Thus $H$ is defined by $M_1 = 1, L_1 = 1$. But, again, this implies $W$ is a strongly geometrically isolated manifold, which is a contradiction. □

**Corollary 6.1.** *For any constant $D > 0$, there are only finitely many Dehn fillings of $W$ whose degrees of trace fields are less than $D$.*

## 7. Final Comments

(1) By Theorem 2.1, for a given $k$-cusped hyperbolic 3-manifold $M$, it easily follows that each pair of equations $M_i = 1, L_i = 1$ produces a maximal anomalous subvariety of its irreducible holonomy variety. Moreover, by Lemma 4.3, if $M$ is a 2-cusped with rationally independent cusp shapes, these are the only maximal anomalous subvarieties of the holonomy variety. Initially, we had thought that Lemma 4.3 would be true for any 2-cusped manifold, but soon realized that it is not true. For instance, if a given 2-manifold has an isometry which sends a cusp to the other (e.g. $W$ in Section 6), then there's a symmetry between $u_1$ and $u_2$ (and between $v_1$ and $v_2$ as well), so we can check $M_1 = M_2, L_1 = L_2$ and $M_1 M_2 = 1, L_1 L_2 = 1$ give other maximal anomalous subvarieties containing $(1, 1, 1, 1)$. As a result, to prove the conjectures in Section 1, it seems that we first need to understand the anomalous subvarieties of the given irreducible holonomy variety.

(2) If a $k$-cusped hyperbolic 3-manifold has a single cusp which is strongly geometrically isolated from the rest, then the irreducible holonomy variety of it is covered by its maximal anomalous subvarieties, meaning that the irreducible holonomy variety is not simple. But, except for this, we don't know of any other example having this property (i.e. non-simple holonomy variety). So, using this fact, we formulate a somewhat stronger conjecture than **Conjecture 1** as follows:

**Conjecture 3.** *Let $X$ be a $k$-cusped hyperbolic 3-manifold. If the deformation variety of $X$ is not simple, then $M$ has a single cusp which is strongly geometrically isolated from the rest.*

Proving or disproving this conjecture would be very interesting and the first step toward **Conjecture 2**.

In general, with the help of Theorem 5.2, to prove the conjectures completely, it is enough to study the unlikely intersections between $\text{Def}(M)$ and linear planes. For this, we need to further investigate properties of higher coefficients of the potential function of a given manifold. However, to the best of our knowledge, we barely have any information about them except for the fact that its first order coefficients (i.e. cusp shapes) are nonreal complex numbers. So we look forward to future research into this direction.

## 8. Proof of Theorem 5.2

In this section, we provide a proof of Theorem 5.2. First set
$$(G_m^{n+t})_*^{[s]} = \bigcup H(\overline{\mathbb{Q}})$$
where the union runs over all algebraic subgroups $H \subset G_m^{n+t}$ of codimension at least $s$ defined by restricting the first $n$-coordinates (or equivalently $\{1\}^n \times G_m^t \subset H$, see Definition 5.1). Then we can restate Theorem 5.2 in the same form as Corollary 1 in [5]:

**Theorem 8.1.** *Let $X \subset G_m^{n+t}$ be an $s$-dimensional variety $(s \leq n)$. Then the height is bounded from above on $X^{re-oa} \cap (G_m^{n+t})_*^{[s]}$.*

To prove Theorem 8.1 we follow exactly the same steps as in [5]. Throughout the proof, $\text{Mat}^*_{s(n+t)}$ means a set of all matrices whose last $t$ columns are all equal to zero (occasionally if we refer to $\text{Mat}_{s(n+t)}$, it simply means the usual set of matrices). Since we only deal with algebraic subgroups (and their cosets) defined by restricting the first $n$-coordinates, we ignore the last $t$-columns. Of course, this plays the same role as "$\text{Mat}_{sn}$" in [5]. Sometimes we denote $(G_m^{n+t})_*^{[s]}$ as $G_*^{[s]}$ for simplicity. All other notation is exactly the same as in [5] unless otherwise stated.

**Lemma 8.2.** *Let $Q > 1$ be a real number and let $\varphi_0 \in \mathrm{Mat}^*_{s(n+t)}(\mathbb{R})$, there exist $q \in \mathbb{Z}$ and $\varphi \in \mathrm{Mat}^*_{s(n+t)}(\mathbb{Z})$ such that*

$$1 \leq q \leq Q \text{ and } |q\varphi_0 - \varphi| \leq \frac{\sqrt{sn}}{Q^{1/(sn)}}.$$

*Proof.* This is essentially the same statement as Lemma 3 in [5]. $\square$

We define $\mathcal{K}^*_{s(n+t)} \subset \mathrm{Mat}^*_{s(n+t)}(\mathbb{R})$ (which corresponds to $\mathcal{K}_{sn}$ in [5]) to be the compact set of all matrices whose rows are orthonormal. All elements of $\mathcal{K}^*_{s(n+t)}$ have rank $s$.

**Lemma 8.3.** *Suppose $W \subset \mathrm{Mat}^*_{s(n+t)}(\mathbb{R})$ is an open neighborhood of $\mathcal{K}^*_{s(n+t)}$. Then there is $Q_0 \geq 1$ (which may depend on $W$) with the following property. For $Q > Q_0$ a real number and $\varphi_0 \in \mathrm{Mat}^*_{s(n+t)}(\mathbb{R})$ with rank $s$, there exist $q \in \mathbb{Z}, \varphi \in \mathrm{Mat}^*_{s(n+t)}(\mathbb{Z})$, and $\theta \in \mathrm{Mat}_s(\mathbb{Q})$ such that*

$$(8.1) \quad 1 \leq q \leq Q, \ \frac{\varphi}{q} \in W, \ |q\theta\varphi_0 - \varphi| \leq \frac{\sqrt{sn}}{Q^{1/(sn)}}, \text{ and } |\varphi| \leq (s+1)q.$$

*Proof.* This is essentially the same statement as Lemma 4 in [5]. $\square$

**Lemma 8.4.** *Let $\varphi : G_m^{n+t} \to G_m^s$ and $p \in G_m^{n+t}(\overline{\mathbb{Q}})$, then*

$$h(\varphi(p)) \leq \sqrt{s(n+t)}|\varphi|h(p).$$

*Proof.* This immediately follows from Lemma 5 in [5]. $\square$

**Lemma 8.5.** *Suppose $W \subset \mathrm{Mat}^*_{s(n+t)}(\mathbb{R})$ is an open neighborhood of $\mathcal{K}^*_{s(n+t)}$. Let $Q_0$ be the constant from Lemma 8.3 and let $Q > Q_0$ be a real number. If $p \in (G_m^{n+t})^{[s]}_*$ then there exist $q \in \mathbb{Z}$ and $\varphi \in \mathrm{Mat}^*_{s(n+t)}(\mathbb{Z})$ such that*

$$1 \leq q \leq Q, \ \frac{\varphi}{q} \in W, \ h(\varphi(p)) \leq \frac{sn}{Q^{1/(sn)}}h(p), \text{ and } |\varphi| \leq (s+1)q.$$

*Proof.* By Lemma 6 in [5], it follows that there exist $q \in \mathbb{Z}$ and $\varphi \in \mathrm{Mat}^*_{s(n+t)}(\mathbb{Z})$ such that

$$1 \leq q \leq Q, \ \frac{\varphi}{q} \in W, \ h(\varphi(p)) \leq \frac{sn}{Q^{1/(sn)}}h(\mathrm{Pr}_n(p)), \text{ and } |\varphi| \leq (s+1)q$$

where $\mathrm{Pr}_n : G_m^{n+t} \longrightarrow G_m^n$ is the projection onto the first $n$-coordinates. Since $h(\mathrm{Pr}_n(p)) \leq h(p)$, we get the desired result. $\square$

From now on, $Y \subset G_m^{n+t}$ denotes an irreducible closed subvariety of $X$ having dimension $r \geq 1$.

The following explanation is given after Lemma 7 in Section 6 in [5], but we repeat it here to make our arguments easy to follow. Note that Habegger denotes our $Y$ as $X$ in his paper. But since we already used $X$ in the statement of Theorem 8.1, we use $Y$ to avoid confusion.

Let $\exp : \mathbb{C}^{n+t} \longrightarrow G_m^{n+t}(\mathbb{C})$ denote the $(n+t)$-fold product of the usual exponential map. It is a locally biholomorphic map between two complex manifolds and as such open. We further assume that $1 := (1)^{n+t}$, the unit element of $G_m^{n+t}$, is a non-singular point of $Y$. Now some open neighborhood $U \subset Y(\mathbb{C})$ of $1$ is an $r$-dimensional complex manifold. After replacing $U$ by a smaller open set we may assume that there is a complex manifold $M \subset \mathbb{C}^{n+t}$ of dimension $r$ containing $0$ such that

$$\exp|_M : M \to U$$

is biholomorphic.

We consider $\varphi \in \mathrm{Mat}_{r(n+t)}(\mathbb{C})$ as a linear map $\mathbb{C}^{n+t} \to \mathbb{C}^r$. Its restriction $\varphi|_M$ is a holomorphic map between two $r$-dimensional complex manifolds. In particular, for each $z \in M$ we have a $\mathbb{C}$-linear differential map

$$d_z(\varphi|_M) : T_z M \to T_{\varphi(z)} \mathbb{C}^r = \mathbb{C}^r$$

between the respective tangent spaces.

**Proposition 8.6.** *Suppose $\varphi : G_m^{n+t} \to G_m^r$ is a nontrivial homomorphism of algebraic subgroups. There exist a dense Zariski open subset $U \subset Y$ and a constant $C_7$ such that*

$$(8.2) \qquad h(\varphi(p)) \geq \frac{r}{2C_1} |\varphi| \frac{\Delta_Y(\varphi)}{|\varphi|^r} h(p) - C_7$$

*for all $p \in U(\overline{\mathbb{Q}})$ where $C_1 = (4(n+t))^r \deg(Y)$.*

*Proof.* This immediately follows from Proposition 1 in [5]. □

**Lemma 8.7.** *Let $\varphi_0 \in \mathrm{Mat}^*_{r(n+t)}(\mathbb{C})$ be such that $d_{z_0}(\varphi_0|_M)$ is an isomorphism of $\mathbb{C}$-vector spaces for some $z_0 \in M$. Then there exist $C_8 > 0$ and an open neighborhood $W \subset \mathrm{Mat}^*_{r(n+t)}(\mathbb{R})$ of $\varphi_0$ such that*

$$\Delta_Y(\varphi) \geq C_8$$

*for all $\varphi \in W \cap \mathrm{Mat}^*_{r(n+t)}(\mathbb{Q})$.*

*Proof.* This easily follows from Lemma 8 in [5]. □

**Lemma 8.8.** *Let $\mathcal{K} \subset \mathrm{Mat}^*_{s(n+t)}(\mathbb{R})$ be compact. One of the following cases holds:*

(1) *There exists $\varphi_0 \in \mathcal{K}$ such that for all $z \in M$ the differential*

$$d_z(\varphi_0|_M) : T_z M \to \mathbb{C}^s$$

*is not injective.*

(2) *There exists $C_9 > 0$ and an open neighborhood $W \subset \mathrm{Mat}^*_{s(n+t)}(\mathbb{R})$ of $\mathcal{K}$ such that for each $\varphi \in W \cap \mathrm{Mat}^*_{s(n+t)}(\mathbb{Q})$ there is $\pi \in \prod_{rs}$ with $\Delta_Y(\pi\varphi) \geq C_9$.*

*Proof.* We will assume that case (i) does not hold and will show that case (ii) does.

Let $\varphi_0 \in \mathcal{K}$. There exist $\pi \in \prod_{rs}$ and $z \in M$ such that $d_z(\pi\varphi_0|M)$ is injective and hence an isomorphism of $\mathbb{C}$-vector spaces. By lemma 8.7 we may find an open neighborhood of $\pi\varphi_0$ in $\mathrm{Mat}^*_{r(n+t)}$ with the stated properties. It follows that we may find $W_{\varphi_0}$, an open neighborhood of $\varphi_0$ in $\mathrm{Mat}^*_{s(n+t)}$, and $C_{\varphi_0}$ with $\Delta(\pi\varphi) \geq C_{\varphi_0}$ for all $\varphi \in W_{\phi_0} \cap \mathrm{Mat}^*_{s(n+t)}(\mathbb{Q})$.

The open cover $\bigcup_{\varphi_0 \in \mathcal{K}} W_{\varphi_0}$ contains $\mathcal{K}$. Since $\mathcal{K}$ is compact we may pass to a finite subcover and conclude that there exist $C_9 > 0$ and an open subset $W$ of $\mathrm{Mat}^*_{s(n+t)}(\mathbb{R})$ containing $\mathcal{K}$ such that for each $\varphi \in W \cap \mathrm{Mat}^*_{s(n+t)}(\mathbb{Q})$ there is $\pi \in \prod_{rs}$ with $\Delta(\pi\varphi) \geq C_9$. □

The following is Proposition 2 in [5].

**Proposition 8.9.** *Let $\mathcal{K} \subset \mathrm{Mat}_{sn}(\mathbb{R})$ be compact and such that all its elements have rank $s$. One of the following cases holds:*

(1) *There exists an algebraic subgroup $H \subset G_m^n$ such that*
$$\dim_p Y \cap pH \geq \max\{1, s + \dim H - n + 1\}$$
*for all $p \in Y(\mathbb{C})$.*

(2) *There exists $C_{10} > 0$ and an open neighborhood $W \subset \mathrm{Mat}_{sn}(\mathbb{R})$ of $\mathcal{K}$ such that for each $\varphi \in W \cap \mathrm{Mat}_{sn}(\mathbb{Q})$ there is $\pi \in \prod_{rs}$ with $\Delta_Y(\pi\varphi) \geq C_{10}$.*

Using the above proposition, we prove an analogous version which we need for Theorem 8.1. This is the key fact that makes it possible to generalize the original Bounded Height Conjecture.

**Proposition 8.10.** *Let $Y \subset X$ be an $r$-dimensional variety such that $\dim Y = \dim \overline{\mathrm{Pr}_n(Y)}$ where $\mathrm{Pr}_n$ is the projection map from $G_m^{n+t}$ to $G_m^n$ (the first $n$-coordinates) and $\mathcal{K} \subset \mathrm{Mat}^*_{s(n+t)}(\mathbb{R})$ be a compact set such that all its elements have rank $s$, then one of the following cases holds:*

(1) *There exists an algebraic subgroup $H \subset G_m^{n+t}$ defined by restricting the first $n$-variables and satisfying*
$$\dim_p Y \cap pH \geq \max\{1, s + \dim H - (n+t) + 1\}$$
*for all $p \in Y(\mathbb{C})$.*

(2) *There exists $C_{10} > 0$ and an open neighborhood $W \subset \mathrm{Mat}^*_{s(n+t)}(\mathbb{R})$ of $\mathcal{K}$ such that for each $\varphi \in W \cap \mathrm{Mat}^*_{s(n+t)}(\mathbb{Q})$ there is $\pi \in \prod_{rs}$ with $\Delta_Y(\pi\varphi) \geq C_{10}$.*

*Proof.* Since the function $\Delta$ is invariant under translation of $Y$, we assume that 1 is a non-singular point of $Y$ (as we previously assumed in the explanation before Proposition 8.6). Moreover, using the assumption $\dim Y = \dim \overline{\mathrm{Pr}_n(Y)}$, we further assume that 1 is a smooth point such that $\mathrm{Pr}_n(1)$ is a smooth point of $\overline{\mathrm{Pr}_n(Y)}$ as well.

If case (2) of Lemma 8.8 holds, then clearly case (2) of this proposition holds. Hence we may assume that we are in case (1) of Lemma 8.8; we will show that case (1) of this proposition holds.

Suppose that there exists $\varphi_0 \in \mathcal{K}$ such that for all $z \in M$ the differential $d_z(\varphi_0|_M)$ fails to be injective. Let $M' \subset \overline{\mathrm{Pr}_n(Y)}$ be a small neighborhood of $\mathrm{Pr}_n(1)$ satisfying $\mathrm{Pr}_n(M) = M'$. If $\varphi'_0 \in \mathrm{Mat}_{sn}(\mathbb{R})$ is the map induced from $\varphi_0$ by removing from the last $t$ columns, then it is easy to check that the differential $d_{z'}(\varphi'_0|_{M'})$ also fails to be injective for all $z' \in M'$. So, by Proposition 8.9, there exists an algebraic subgroup $H \subset G_m^n$ such that $\overline{\mathrm{Pr}_n(Y)}$ is equal to

$$\left\{ p \in \overline{\mathrm{Pr}_n(Y)} \mid \dim_p \overline{\mathrm{Pr}_n(Y)} \cap pH \geq \max\{1, s + \dim H - n + 1\} \right\}.$$

Thinking of $H$ as an algebraic subgroup in $G_m^{n+t}$, the above fact implies that the closed subvariety

$$\{p \in Y \mid \dim_p Y \cap pH \geq \max\{1, s + \dim H - (n+t) + 1\}\}$$

of $Y$ contains $M$, which means

$$Y = \{p \in Y \mid \dim_p Y \cap pH \geq \max\{1, s + \dim H - (n+t) + 1\}\}.$$

This completes the proof. $\square$

The following lemma enables us to apply the above proposition to any closed irreducible subvariety $Y \subset X$ such that $Y \cap X^{re-oa} \neq \emptyset$.

**Lemma 8.11.** *If $Y \subset X$ is an irreducible closed subvariety such that $Y \cap X^{re-oa} \neq \emptyset$, then $\dim Y = \dim \overline{\mathrm{Pr}_n(Y)}$.*

*Proof.* Let $\dim Y - \dim \overline{\mathrm{Pr}_n(Y)} = k \ (\geq 1)$. Then

$$Y = \{p \in Y(\mathbb{C}) : \dim_p Y \cap p(\{1\}^n \times G_m^t) \geq k\}$$

(for example, see Theorem 3.13 in [7]). Since the set $\{p \in Y(\mathbb{C}) : \dim_p Y \cap p(\{1\}^n \times G_m^t) \geq k\}$ is in the complement of $X^{re-oa}$, we get the desired result. $\square$

In the proof of the next lemma, we simply copy the proof of Lemma 11 in [5] except for adjusting the constants given in (8.4). Before proceeding we make the following easy observation: say $\epsilon$ is an small number satisfying $0 < \epsilon \leq \frac{1}{2(n+t)}$ with $p \in \mathcal{C}(G_*^{[s]}, \epsilon)$, so there is $a \in G_*^{[s]}$ and $b \in G_m^{n+t}(\overline{\mathbb{Q}})$ with $h(b) \leq \epsilon(1 + h(a))$. Then $h(a) = h(pb^{-1}) \leq h(p) + h(b^{-1}) \leq h(p) + (n+t)h(b) \leq h(p) + (1 + h(a))/2$ by the elementary properties of height. We easily deduce

(8.3) $$h(a) \leq 1 + 2h(p), \quad h(b) \leq 2\epsilon(1 + h(p)).$$

**Lemma 8.12.** *Let $Y \subset X$ be an irreducible closed subvariety of positive dimension. If $Y \cap X^{re-oa} \neq \emptyset$, there exists $\epsilon > 0$ and $U \subset Y$ which is Zariski open and dense such that the height is bounded on $U(\overline{\mathbb{Q}}) \cap \mathcal{C}(G_*^{[s]}, \epsilon)$.*

*Proof.* Since $Y \cap X^{re-oa} \neq \emptyset$, by Lemma 8.11 it satisfies $\dim Y = \dim \overline{\mathrm{Pr}_n(Y)}$. So we can apply Proposition 8.10, and the same condition (i.e. $Y \cap X^{re-oa} \neq \emptyset$) means that we are in case (2) of Proposition 8.10 applied $\mathcal{K} = \mathcal{K}^*_{s(n+t)}$.

Therefore there exist an open set $W \subset \mathrm{Mat}^*_{s(n+t)}(\mathbb{R})$ containing $\mathcal{K}^*_{s(n+t)}$ and $C_{10} > 0$ such that for each $\varphi \in W \cap \mathrm{Mat}^*_{s(n+t)}(\mathbb{Q})$ there is $\pi \in \prod_{\dim Y, s}$ with $\Delta_Y(\pi\varphi) \geq C_{10}$.

We suppose $Q_0$ is as in Lemma 8.5 and that $Q > Q_0$ is a fixed parameter which depends only on $X$ and $Y$. We will see later how to choose $Q$ properly.

Let $\Theta$ denote the set of all matrices $\varphi \in \mathrm{Mat}^*_{s(n+t)}(\mathbb{Z})$ such that there exists an integer $q$ with $1 \leq q \leq Q, \varphi/q \in W$, and $|\varphi| \leq (s+1)q$ (cf. Lemma 8.5). Clearly, $\Theta$ is a finite set.

For each $\varphi \in \Theta$ there is a $\pi \in \prod_{\dim Y, s}$ such that $\Delta_Y(\varphi'/q) \geq C_{10}$ where $\varphi' = \pi\varphi$. In particular, $\varphi' \neq 0$ since $C_{10} > 0$. By homogeneity we have
$$\Delta_Y(\varphi') = q^{\dim Y} \Delta_Y(\varphi'/q) \geq C_{10} q^{\dim Y}.$$
Now $\varphi' \neq 0$ implies $|\varphi'| \geq 1$ so we obtain the following lower bound for the factor in front of $h(p)$ in (8.2)
$$C_{11}|\varphi'|\frac{\Delta_Y(\varphi')}{|\varphi'|^{\dim Y}} \geq C_{10}C_{11}|\varphi'|\frac{q^{\dim Y}}{|\varphi'|^{\dim Y}} \geq C_{10}C_{11}\frac{q^{\dim Y}}{|\varphi'|^{\dim Y}}$$
with
$$C_{11} = \frac{\dim Y}{2(4(n+t))^{\dim Y} \deg(Y)} > 0.$$
Now $|\varphi'| \leq |\varphi| \leq (s+1)q$, so
$$C_{11}|\varphi'|\frac{\Delta_Y(\varphi')}{|\varphi'|^{\dim Y}} \geq \frac{C_{10}C_{11}}{(s+1)^{\dim Y}}.$$
We denote this last quantity by $C_{12}$; it is positive and independent of $Q$ and $\varphi$.

We fix
$$Q = \max\left\{Q_0 + 1, (8s(n+t)C_{12}^{-1})^{sn}\right\} > Q_0,$$
(8.4)
$$\epsilon = \min\left\{\frac{1}{2(n+t)}, \frac{\sqrt{s(n+t)}}{s+1}\frac{1}{Q^{1+1/(sn)}}\right\} \in \left(0, \frac{1}{2(n+t)}\right].$$

Let $U_\varphi$ be the dense Zariski open subset of $Y$ supplied by Proposition 8.6 applied to $\varphi$. The intersection
$$U = \bigcap_{\varphi \in \Theta} U_\varphi$$
is a dense Zariski open subset of $Y$ since $\Theta$ is finite. We deduce that
(8.5)
$$h(\varphi'(p)) \geq C_{12}h(p) - C(Q)$$
for all $p \in U(\overline{\mathbb{Q}})$ and all $\varphi \in \Theta$; here $C(Q)$ depends neither on $p$ nor on $\varphi$ (but possibly on $Q$).

Now let us assume that $p \in U(\overline{\mathbb{Q}}) \cap \mathcal{C}(G_*^{[s]}, \epsilon)$. That is, there are $a \in G_*^{[s]}$ and $b \in G_m^{n+t}(\overline{\mathbb{Q}})$ with $p = ab$ and $h(b) \leq \epsilon(1 + h(a))$.

By Lemma 8.5 there exists $\varphi \in \Theta$ with $h(\varphi(a)) \leq snQ^{-1/(sn)}h(a)$ and so
(8.6)
$$h(\varphi(a)) \leq 2snQ^{-1/(sn)}(1 + h(p))$$

by (8.3).

We apply Lemma 8.4 in order to bound $h(\varphi(b)) \leq \sqrt{s(n+t)}|\varphi|h(b)$. Now (8.3) gives $h(\varphi(b)) \leq 2\epsilon\sqrt{s(n+t)}|\varphi|(1+h(p))$. But $|\varphi| \leq (s+1)q \leq (s+1)Q$, so

$$(8.7) \qquad h(\varphi(b)) \leq 2\epsilon\sqrt{s(n+t)}(s+1)Q(1+h(p)).$$

Using (8.6), (8.7), and elementary properties of height give

$$h(\varphi(p)) = h(\varphi(ab)) \leq h(\varphi(a)) + h(\varphi(b))$$
$$\leq \left(2snQ^{-1/(sn)} + 2\epsilon\sqrt{s(n+t)}(s+1)Q\right)(1+h(p)).$$

The choice of $\epsilon$ made in (8.4) implies $h(\varphi(p)) \leq 4s(n+t)Q^{-1/(sn)}(1+h(p))$ and the choice of $Q$ gives $h(\varphi(p)) \leq C_{12}(1+h(p))/2$. Furthermore, we have $h(\varphi'(p)) \leq h(\varphi(p))$, hence

$$(8.8) \qquad h(\varphi'(p)) \leq \frac{C_{12}}{2}(1+h(p)).$$

If we compare (8.5) and (8.8) we immediately get the desired $h(p) \leq 1 + 2C_{12}^{-1}C(Q)$. $\square$

For brevity we set $\Sigma = X^{re-oa} \subset X(\overline{\mathbb{Q}})$. If $X^{re-oa} \neq \emptyset$, then Lemma 8.12 applied with $X = Y$ shows that there exists a dense Zariski open subset $U \subset X$ such that $U(\overline{\mathbb{Q}}) \cap G_*^{[s]}$ has bounded height. This is already close to Theorem 8.1 and the following simple descent argument shows how to deal with the points in $(\Sigma \setminus U(\overline{\mathbb{Q}})) \cap G_*^{[s]}$:

**Lemma 8.13.** *Suppose that there is a proper subset $S \subsetneq \Sigma$ and an $\epsilon > 0$ such that the height is bounded from above on $S \cap \mathcal{C}(G_*^{[s]}, \epsilon)$. There exists a subset $S' \subset \Sigma$ containing $S$ with $\overline{\Sigma \setminus S'} \subsetneq \overline{\Sigma \setminus S}$ and an $\epsilon' > 0$ such that the height is bounded from above on $S' \cap \mathcal{C}(G_*^{[s]}, \epsilon')$.*

*Proof.* This easily follows by copying the proof of Lemma 12 in [5]. $\square$

*Proof of Theorem 8.1.* This also easily follows from the proof given in [5]. $\square$

## 9. Acknowledgement

I would like to thank to my advisor Nathan Dunfield for all his support and guidance throughout this project. All the conversations with him helped me to improve every result and proof of the paper. I'm also very much indebted to Professor Umberto Zannier for introducing and explaining Theorem 1.2 and Theorem 3.9 as well as pointing out an error in the earlier version of the paper.

Department of Mathematics
University of Illinois at Urbana-Champaign
1409 W. Green Street, Urbana, IL 61801

*Email Address*: jeon14@illinois.edu